\documentclass{article}
\usepackage{latexsym}
\usepackage{amsmath}
\usepackage{amssymb}

\newcommand{\Q}{\mathbb{Q}}

\newcommand{\aQ}{\overline{\mathbb{Q}}}
\newcommand{\Z}{\mathbb{Z}}
\newcommand{\Zp}{\mathbb{Z}_{p}}

\newcommand{\aQp}{\overline{\mathbb{Q}_{p}}}
\newcommand{\Qp}{\mathbb{Q}_{p}}
\newcommand{\GL}{\mathrm{GL}}
\newcommand{\Gm}{\mathbb{G}_{\mathrm{m}}}

\newtheorem{prop}{Proposition}
\newtheorem{theo}{Th\'eor\`eme}
\newtheorem{cor}{Corollaire}
\newtheorem{lemme}{Lemme}

\title{Sur les repr\'esentations $p$-adiques g\'eom\'etriques 
de conducteur $1$ et de dimension $2$ de $G_{\Q}$.}

\author{J.-P. ~\textsc{Wintenberger}}

\begin{document}

\maketitle

\it R\'esum\'e. \rm En bas poids et sous des hypoth\`eses
convenables,
nous prouvons qu'il n'existe pas
de repr\'esentation $p$-adique g\'eom\'etrique irr\'eductible
de conducteur $1$ et de dimension $2$ de $G_{\Q}$, conform\'ement
\`a une conjecture de Fontaine et Mazur.

\it Summary. \rm 
We prove that there is no geometric $p$-adic 
representation of the Galois group of $\Q$ which is irreducible,
of dimension $2$, of conductor $1$ and low weight, 
according to a conjecture of Fontaine and Mazur.

\section{Introduction et \'enonc\'e du r\'esultat principal.} 

Soient $p$ un nombre premier et  $\Qp$ le corps des nombres
$p$-adiques.
Soit $\aQ$ une cl\^oture alg\'ebrique
de $\Q$.
On note  $G_{\Q}$  le groupe de Galois de $\aQ / \Q$.
Pour $L$ extension de $\Q$ contenue dans $\aQ$,
on note $G_L$ le groupe de Galois de $\aQ / L$.
On d\'esigne par $E$ une extension 
finie de $\Qp$. Une repr\'esentation $p$-adique de $G_{\Q}$
\`a coefficients dans $E$ est un homomorphisme continu $\rho$ de 
$G_{\Q}$  dans $\GL _E (U)$, o\`u $U$ est
un $E$-espace vectoriel de 
dimension finie.

Une telle repr\'esentation est \it g\'eom\'etrique \rm si elle v\'erifie
les deux conditions suivantes (\cite{[FM95]}):

- sa restriction \`a un sous-groupe de d\'ecomposition
$D_{p}$ en $p$  est potentiellement semi-stable 
au sens de la th\'eorie de Fontaine (exp. 8 de \cite{[Per]}),

- il existe un ensemble fini $S$ de nombres premiers tel 
que $\rho$ soit non ramifi\'ee en dehors de $S\cup \{ p\}$.

Une  repr\'esentation g\'eom\'etrique a un conducteur $N(\rho )$
(\cite{[FM95]}). Pour 
$\ell \not= p$, l'action du groupe de d\'ecomposition $D_{\ell}$
sur $U$ d\'efinit une action du groupe de Weil-Deligne $\mathrm{WD}_{\ell}$
sur $U$, donc la composante $N_{\ell}=l^{*}$ en  $\ell$ du conducteur
$N(\rho )$ ($N_{\ell}=1$ si $\ell\notin S\cup \{ p \}$).
On a de m\^eme une action du groupe de Weil-Deligne $\mathrm{WD}_{p}$
sur le module de Dieudonn\'e filtr\'e associ\'e \`a la restriction
de $\rho$ \`a $D_{p}$, ce qui d\'efinit la composante en $p$
de $N(\rho )$. En particulier, $N(\rho )=1$ si et seulement
si $\rho$ est cristalline
en $p$ et non ramifi\'ee hors de $p$.

Soit $\rho : G_{\Q}\rightarrow \GL _E  (U)$ une repr\'esentation
$p$-adique \`a coefficients dans $E$.   
Il n'est pas difficile de voir que, si la restriction
de $\rho$ \`a un sous-groupe de d\'ecomposition $D_p$ en $p$
est de Hodge-Tate, les 
poids de Hodge-Tate de $\rho$, consid\'er\'ee comme
repr\'esentation dans le $\Qp$-espace vectoriel $U$,
ont des multiplicit\'es qui sont multiples de $[E:\Qp ]$
(prop \ref{htconstant}).
On appelle \it poids de Hodge-Tate  \rm
de $\rho$ ces poids, avec les multiplicit\'es 
divis\'ees par $[E:\Qp ]$.

Supposons $\rho$ de dimension $2$. 
Soit $c\in G_{\Q}$ une conjugaison complexe. La 
repr\'esentation $\rho$ est \it  impaire \rm  si $\rho (c)$
a comme valeurs propres $1$ et $-1$, \it i.e. \rm si $\rho (c)$ est  
de d\'eterminant $-1$.

Soit $f=q+\ldots + a_{n} q^n+\ldots$ une forme modulaire primitive 
(parabolique, propre pour les op\'erateurs de Hecke) 
de conducteur $N(f)$ (sur $\Gamma _{1} (N(f) )$).
Notons $E(f )$ son corps de rationalit\'e, \it i.e. \rm le corps
engendr\'e par les coefficients de $f$ et les valeurs du caract\`ere
de $f$. Le corps $E(f )$ est une extension finie de 
$\Q$. 
Soit $\aQp$ une cl\^oture alg\'ebrique de $\Qp$.
On sait associer \`a $f$ et \`a un  plongement  
$i$ de $E (f)$ dans $\aQp$ une repr\'esentation
$p$-adique $\rho _{f,i}: G_{\Q}\rightarrow 
\GL _{2} (E)$, $E$ \'etant  l'adh\'erence
de l'image de $i$. Elle est non ramifi\'ee en $\ell$
pour tout nombre premier $\ell$ qui est distinct de $p$
et ne divise pas $N(f)$ et est  caract\'eris\'ee par :

$$ \mathrm{tr}(\rho (\mathrm{Frob}_{\ell}))= i(a_{\ell}),$$
 
pour ces $\ell$. De plus, $\rho _{f,i}$ est absolument 
irr\'eductible, impaire, g\'eom\'etrique  de conducteur
$N(f)$ et de poids de Hodge-Tate $(0,k-1)$. 

Il est conjectur\'e que, si 
 $\rho :G_{\Q} \rightarrow \GL _{2} (E)$ est une repr\'esentation
absolument irr\'eductible, impaire, g\'eom\'etrique 
de conducteur $N$ et de poids de Hodge-Tate $(0,k-1)$, alors $\rho$ 
est isomorphe \`a $\rho_{f,i}$ pour $f$ forme primitive
de conducteur  $N$ et de poids $k$ et $i$ un plongement 
du corps de rationalit\'e $E( f)$ de $f$  dans $E$
(\cite{[FM95]}).

Une forme primitive de conducteur $1$ a un poids $k$ pair.
Soit $\rho : G_{\Q}
\rightarrow \GL _{2} (E )$ g\'eom\'etrique de conducteur $1$   
et de poids $(0,k-1)$.
Il n'est pas difficile de prouver que  $\mathrm{det}(\rho )=\chi _{p}^{k-1}$,
$\chi _{p} : G_{\Q}\rightarrow \Qp ^{*}$ d\'esignant le caract\`ere
cyclotomique (corollaire \ref{det}). On voit donc que $\rho$ est impaire 
si et seulement si $k$ est pair.
Comme il n'existe pas de forme primitive $f$
de conducteur $1$ et et poids $k< 12$, la conjecture de Fontaine
et Mazur pr\'edit
qu'il n'existe pas de $\rho$ g\'eom\'etrique,
absolument irr\'eductible de conducteur $1$ et de poids 
$(0, k-1)$ pour $k$ pair $< 12$.

Fontaine et Abrashkin ont prouv\'e que, pour des $p$
petits, il n'existe pas de 
repr\'esentation $p$-adique de $G_{\Q}$ qui soit g\'eom\'etrique
de conducteur $1$, de poids de Hodge-Tate petits, irr\'eductible
de dimension $\not=1$ (\cite{[F85]},\cite{[F93]},\cite{[A89]}). 
Taylor a prouv\'e une version potentielle de la conjecture
de Fontaine et Mazur (\cite{[T02]},\cite{[T03]}).
L'objet de cet article est de montrer  que les th\'eor\`emes de Taylor, 
et de Fontaine
et Abrashkin, permettent de prouver dans certains cas
cette inexistence de repr\'esentation
de conducteur $1$ et dimension $2$.      
 
En fait, on a un \'enonc\'e concernant des repr\'esentations
qui ne sont pas n\'ecessairement semi-simples. 
Soit   $H^1_f (\Q , \Q_p (m))$ le
$\Q _p$-espace vectoriel
qui classe les 
extensions de $\Qp$ par $\Qp (m )$ 
qui sont de conducteur $1$.
On sait  que $H^1_f (\Q , \Q_p (m))$
est de dimension $1$ si
$m$ est impair $\geq 3$ et nul sinon
(ceci r\'esulte essentiellement de travaux de  Soul\'e :
remarque du 4.2. de \cite{[BPR]}, lemme 4.3.1. de \cite{[BN]}).
Il existe donc une repr\'esentation $p$-adique $U_{p,3}$ de $G_{\Q}$ 
qui est g\'eom\'etrique de conducteur $1$
et qui est une extension non triviale
d'une repr\'esentation de dimension $1$ avec action de $G_{\Q}$
triviale, par une repr\'esentation de dimension $1$ avec 
action de  $G_{\Q}$ donn\'ee par $\chi _p ^3$. 
Cette repr\'esentation est unique \`a isomorphisme
pr\`es. Nous prouvons :  

\begin{theo}\label{theoprinc} Soient $E$ une extension finie
de $\Qp$, 
$U$ un $E$-espace vectoriel de dimension $2$ et 
$\rho : G_{\Q}\rightarrow
\GL_E (U)$ une repr\'esentation $p$-adique \`a coefficients 
dans $E$.
On suppose que $\rho$ est g\'eom\'etrique 
de conducteur $1$ de poids de Hodge-Tate
 $(0,k-1)$ avec $k=2$ ou $k=4$  et que dans ce deuxi\`eme
cas, que $p$ soit $\geq 7$. Alors, $\rho$ est isomorphe :

- si $k=2$, \`a $E \oplus E (1)$ (l'action $G_{\Q}$ sur $E$ \'etant 
triviale) ;

- si $k=4$, soit \`a $E \oplus E  (3)$, soit \`a $E\otimes_{\Qp } U_{p,3}$.
\end{theo} 

\it Remarque. \rm Le cas $k=1$ est sans int\'er\^et, puisque si $\rho : G_{\Q}
\rightarrow \GL _{2} (E)$ est g\'eom\'etrique
de conducteur $1$ et de poids
$(0,0)$, elle est partout non ramifi\'ee, donc triviale.
Il est possible que l'on puisse de d\'ebarrasser dans le th\'eor\`eme
de l'hypoth\`ese $p\geq 7$ par les m\'ethodes 
de Fontaine et Abrashkin. 

\vskip0.5cm

Le coeur de la preuve est au 2.4. la construction d'une repr\'esentation
$q$-adique, $q$ nombre premier $\not=p$, qui est potentiellement
compatible \`a $\rho$. Cette construction repose sur le th\'eor\`eme
de Taylor et un petit compl\'ement \`a ce th\'eor\`eme (prop. 
\ref{complement}).
On conclut gr\^ace aux th\'eor\`emes de Fontaine et Abrashkin.
Le 2.1. et le 2.2 sont des pr\'eliminaires. Dans le 3, nous donnons
une application aux groupes $p$-divisibles sur $\Z$. Dans le 
4, nous prouvons que des repr\'esentations $p$-adiques g\'eom\'etriques
de conducteur $1$ de bas poids n'existent pas, sous une hypoth\`ese
d'ordinarit\'e en $3$ : nous utilisons un th\'eor\`eme
de Serre (\cite{[Se]} : Oeuvres, vol 3 p. 710 )
au lieu des th\'eor\`emes de Fontaine et Abrashkin
et un th\'eor\`eme de Skinner et Wiles (\cite{[SW99]}).
  
Je voudrais remercier Henri Carayol pour avoir r\'epondu 
\`a mes multiples questions. Je voudrais aussi remercier Laurent 
Berger, 
Christophe Breuil, Mladen Dimitrov, Jean-Marc Fontaine,  
Eknath Ghate, Bernadette Perrin-Riou et
Jacques Tilouine pour des conversations utiles.  
Alors que je terminais la r\'edaction de ce papier,
Christophe Breuil   m'a signal\'e une pr\'epublication
de Luis Dieulefait avec des r\'esultats similaires
aux notres (\cite{[Dieu]}).

\section{Preuve du th\'eor\`eme \ref{theoprinc}.}\label{preuve}

\subsection{Repr\'esentations de Hodge-Tate
de $G_{\Qp}$.}\label{ht}

Soit $\aQp$ une cl\^oture alg\'ebrique de $\Qp$ ;
d\'esignons par $G_{\Qp}$ le groupe de Galois de
$\aQp / \Qp$ et par $C_p$ le compl\'et\'e de $\aQp$.

\subsubsection{Poids des repr\'esentations de Hodge-Tate
de $G_{\Qp}$.}\label{htpoids}

Si  $G'$ est un sous-groupe ouvert de $G_{\Qp}$ et 
$W$ un $C_p$-espace vectoriel de dimension finie 
muni d'une action semi-lin\'eaire et continue de 
$G'$, $W$ admet une d\'ecomposition
de Hodge-Tate si l'on 
a :

$$ W=\oplus_{i\in \Z}( C_p \otimes_{\aQp^{G'}} W(-i)^{G'}),$$ 

$W(-i)$ \'etant $W$ muni de l'action de $G'$ tordue
par $\chi _p ^{-i}$ (\cite{[Ta67]}). Les poids  de Hodge-Tate de $W$ 
sont  alors  les  entiers $i$ tels que $ W(-i)^{G'}\not= (0)$,
compt\'es avec la multiplicit\'e $\mathrm{dim}_{\aQp^{G'}} (W(-i)^{G'})$.

Soit $U$ un $E$-espace vectoriel de dimension finie $d$.
Soit $\rho : G'\rightarrow \GL _E (U)$ une repr\'esentation
$p$-adique \`a coefficients dans $E$. 
On suppose que $\rho$ est de Hodge-Tate
\it i.e.\rm  que $W=C_p \otimes_{\Qp} U$, muni de l'action 
de $G'$ produit
tensoriel, admet une d\'ecomposition de Hodge-Tate.
Notons $(\lambda ,w)\mapsto [\lambda ]w$, pour $\lambda \in E$ et $w\in W$,
l'action de $E$ sur $W$ induite par celle de $E$ sur $U$.
Notons $P(E,\aQp )$ l'ensemble des plongements de 
$E$ dans $\aQp$. Pour $\iota \in P(E,\aQp )$, 
notons $W_{\iota}$ le sous-espace vectoriel de $W$
form\'e des $w\in W$ v\'erifiant $[\lambda ]w=\iota (\lambda )w$. 
On a donc :

$$ W=\oplus_{P(E,\aQp )} W_{\iota}.$$

Un sous-espace vectoriel de $W$ 
stable par $G'$ admet une d\'ecomposition de Hodge-Tate
(\it cf \rm par exemple exp. 3 de \cite{[Per]}).
Il en r\'esulte que
chacun des  $W_{\iota}$
admet une d\'ecomposition de Hodge-Tate. 
On appelle \it type de Hodge-Tate de $\rho$ \rm
l'application de l'ensemble $P(E,\aQp )$
dans les familles de $d$ entiers qui \`a $\iota$
associe les poids de Hodge-Tate de $C_p \otimes_{\iota, E} U$.
La proposition suivante est   bien connue :

\begin{prop}\label{htconstant} 
On suppose $G'=G_{\Qp}$. Alors, le type de Hodge-Tate
de $\rho$ est constant : les poids de Hodge-Tate de 
$C_p \otimes_{\iota, E} U$ ne d\'ependent pas de $\iota$.
\end{prop}

On appelle \it poids de Hodge-Tate \rm de la  
repr\'esentation $p$-adique $\rho$ de $G_{\Qp}$ \`a coefficients dans $E$  
les poids de Hodge-Tate de  $C_p \otimes_{\iota, E} U$, pour un $\iota$
(compt\'es avec leurs multiplicit\'es dans $C_p \otimes_{\iota, E} U$).

\vskip 0.5cm

\it Preuve de la proposition. \rm 
Reprenons les notations ci-dessus.
Soit $\sigma \in G_{\Qp}$. Comme, si $w\in W_{\iota}$ :

$$ [\lambda ](\sigma w)=\sigma ([\lambda ]w)=
\sigma (\iota(\lambda )w)=(\sigma \iota)(\lambda )\sigma (w),$$

on voit que $\sigma$ d\'efinit une bijection de 
$W_{\iota}$ sur $W_{\sigma \iota}$. Soit $G''$
le sous-groupe ouvert de $G_{\Qp}$ qui fixe le sous-corps
de $\aQp$ engendr\'e par les $\iota (E)$. Le groupe $G''$ 
agit sur $W_{\iota}$ et  $W_{\sigma \iota}$.
Notons $\mathrm{int}(\sigma^{-1} )(W_{\iota})$ le $C_p$-espace
vectoriel $W_{\iota}$, muni de l'action de $G''$ d\'efinie
par $(\tau,w)\mapsto \sigma^{-1} \tau \sigma (w)$. 
Alors, $\sigma$ induit un isomorphime des $C_p$-espaces vectoriels
munis des actions semi-lin\'eaires de $G''$  :

$$ C_p \otimes _{\sigma,C_p} \mathrm{int}(\sigma^{-1} )(W_{\iota})
\simeq W_{\sigma \iota}.$$ 

Comme $C_p \otimes _{\sigma,C_p}
\mathrm{int}(\sigma^{-1} )(W_{\iota}) $
a les m\^emes poids de Hodge-Tate que $W_{\iota}$,
la proposition en r\'esulte.$\Box$

\subsubsection{Le cas des caract\`eres.}\label{caracteres}

\begin{prop}\label{ab} Soit $\eta : G_{\Q}\rightarrow E^*$
un caract\`ere continu de $G_{\Q}$. On suppose que 
$\eta$ est non ramifi\'e en dehors de $p$ et que 
la repr\'esentation $p$-adique $E(\eta)$ de dimension $1$
qu'il d\'efinit est cristalline en $p$. Alors, il existe
un entier $i\in \Z$ tel que $\eta$ soit le 
compos\'e de $\chi_p ^i : G_{\Q}\rightarrow \Qp ^*$
avec l'inclusion $\Qp ^* \hookrightarrow E^*$.
\end{prop} 

\it D\'emonstration.\rm Soit $i$ le poids de Hodge-Tate de 
$E(\eta)$ (\it cf \rm proposition pr\'ec\'edente).
La repr\'esentation  $E(\eta \chi_p ^{-i})$
est non ramifi\'ee partout, donc triviale.$\Box$

\begin{cor}\label{det} Soit $\rho : G_{\Q}\rightarrow 
\GL _2 (E)$ une repr\'esentation g\'eom\'etrique
de conducteur $1$ et de poids de Hodge-Tate $(0,k-1)$.
On a $\mathrm{det}(\rho )=\chi_p ^{k-1}$.\end{cor}

\subsection{Le cas des repr\'esentations potentiellement 
ab\'eliennes.}\label{potab}

\begin{prop}\label{reppotab} Soit $\rho : G_{\Q}
\rightarrow \GL  _{d} (E)$ une repr\'esentation 
$p$-adique qui est irr\'eductible, g\'eom\'etrique
de conducteur $1$ et qui est potentiellement ab\'elienne :
il existe $L\subset \aQ$ extension finie de $\Q$
telle que la restriction de $\rho$ au groupe de Galois $G_{L}$
ait une image ab\'elienne. Alors, $d=1$ et
$\rho$ est isomorphe \`a $E(j)$, pour un entier $j$. \end{prop}

Prouvons tout d'abord le lemme :

\begin{lemme}\label{imconnexe} Soit $\eta : G_{\Q}
\rightarrow \GL  _{d'} (\Qp)$ une repr\'esentation
$p$-adique g\'eom\'etrique de conducteur $1$. Alors,
l'adh\'erence de Zariski de l'image de $\eta$ est connexe.
\end{lemme}

\it Preuve du lemme. \rm Comme $\Q$ n'a pas d'extension
partout non ramifi\'ee, $\eta ( G_{\Q} )$  est engendr\'ee
par les images des sous-groupes d'inertie au dessus de $p$.
Comme $\rho$ est cristalline en $p$, ces images ont  
une adh\'erence de Zariski connexe
(prop. 3.8.4. de  \cite{[F78]}). Le lemme en r\'esulte.

\it Remarque. \rm Le lemme vaut si $\eta$ est seulement 
suppos\'ee non ramifi\'ee hors de $p$ et semi-stable en $p$.

\vskip0.5cm
Prouvons la proposition. Soit $\eta$ la repr\'esentation
$\rho$ vue comme repr\'esentation dans $\GL_{d[E:\Qp ]} (\Qp )$.
Soient $H$ l'adh\'erence de Zariski de l'image de $\eta$ et 
$H^{0}$ sa composante neutre. Comme $\eta$ est semi-simple et potentiellement
ab\'elienne, $H^{0}$ est un tore. Le lemme dit que $H=H^0$, donc 
$\eta$ est ab\'elienne. Elle est localement alg\'ebrique au sens de 
Serre (3.1. de \cite{[SAB]}). Comme elle est de conducteur
$1$, elle provient d'une repr\'esentation du 
groupe de type multiplicatif $S_{\Q,1}$,
relatif \`a $\Q$ et de conducteur $1$ (2 de \cite{[SAB]}).
Le  $S_{\Q,1}$ est r\'eduit \`a $\Gm$,
la repr\'esentation $G_{\Q}\rightarrow \Qp ^{*}$ \'etant donn\'ee 
par le caract\`ere cyclotomique. Comme $\eta$ est obtenue en composant 
cette repr\'esentation avec une repr\'esentation de $\Gm$, on voit 
que $\eta$ est somme de $\Qp  (j)$. La proposition en r\'esulte.$\Box$

\subsection{Le th\'eor\`eme de Taylor.}\label{taylor}

\subsubsection{Enonc\'e.}

Dans \cite{[T02]} et \cite{[T03]}, Taylor prouve :

\begin{theo}\label{taylortheo} On suppose  $p>3$. Soit $k$ un 
entier.
Soit $\rho : G_{\Q} \rightarrow \GL_2 (E)$ une repr\'esentation
$p$-adique absolument 
irr\'eductible , impaire, non ramifi\'ee en dehors d'un ensemble
fini $S$ de premiers, et cristalline en $p$ de poids de Hodge-Tate
$(0,k-1)$. On suppose que $2\leq k \leq (p+1)/2$.
Alors, il existe une extension finie $F\subset \aQ$ de $\Q$ totalement 
r\'eelle, galoisienne, non ramifi\'ee en $p$, et une 
repr\'esentation alg\'ebrique cuspidale  r\'eguli\`ere $\pi$ de 
$\GL_2 (\mathbb{A}_F )$ de poids $k$, non ramifi\'ee en les 
premiers de $F$ au dessus de $p$, 
et un plongement $i$ du corps
de rationalit\'e $E(\pi )$ de $\pi$ dans $E$ tels que 
la restriction de $\rho$ \`a $F$ soit isomorphe
\`a $\rho_{\pi,i}$.
\end{theo}

\it Remarque. \rm R Taylor prouve un \'enonc\'e plus g\'en\'eral
o\`u l'hypoth\`ese $2\leq k \leq (p+1)/2$ est remplac\'ee par
$2\leq k \leq p-1$ avec des hypoth\`eses d'irr\'eductibilit\'e 
de la restriction de la r\'eduction de $\rho$ au corps 
quadratique non ramifi\'e en dehors de $p$ (th. 6.1.
de \cite{[T03]}). 

\subsubsection{Compl\'ement.}

Nous avons besoin du compl\'ement suivant :

\begin{prop}\label{complement} Soient $l_1,\ldots, l_r$ 
des premiers $\not= p$  et  n'appartenant pas \`a $S$.
Alors, on peut supposer $F$  non ramifi\'e  au dessus de  
$l_1,\ldots, l_r$.
\end{prop}

\it Preuve. \rm 
On reprend la preuve de  Taylor.

Soit $\overline{\rho}$ une r\'eduction de $\rho$.
Si  $\overline{\rho}$ n'est pas irr\'eductible ou que 
$\overline{\rho}$ est irr\'eductible \`a image r\'esoluble 
et que la restriction de $\overline{\rho}$
au groupe d'inertie $I_p$ en $p$ a des caract\`eres
de niveau $1$,
le th\'eor\`eme
de Taylor est vrai avec $F=\Q$ d'apr\`es Skinner et 
Wiles puisque $\overline{\rho}$ est $D_p$
distingu\'ee et $\rho$ ordinaire (\cite{[SW99]},\cite{[SW01]},
\cite{[Sk03]}).

Supposons que $\overline{\rho}$ a une image qui n'est pas 
r\'esoluble et que la restriction de $\overline{\rho}$ \`a $I_p$
a des caract\`eres de niveau $1$. On est alors dans le cas de 
\cite{[T02]}. Dans \cite{[T02]}, la repr\'esentation $\rho$
est une repr\'esentation $\ell$-adique (et pas $p$-adique), et 
pour la suite de cette preuve, nous reprenons la notation
de Taylor. La preuve de Taylor utilise un premier $p$
(p. 131) choisi par application du th\'eor\`eme de 
Chebotarev et en dehors d'un ensemble fini de premiers ;
on peut supposer $p$ distinct des $l_i$.

Dans le lemme 1.1. de \cite{[T02]}, si on suppose $L$ non ramifi\'e 
au-dessus des $l_i$, $S$ ne contenant pas de premiers
de $K$ de caract\'eristique r\'esiduelle l'un des $\ell_i$,
et le caract\`ere $\phi$ non ramifi\'e en les premiers de $K$
au dessus des $l_i$, on peut choisir le caract\`ere 
$\Psi$ non ramifi\'e en les premiers de $L$
au dessus des $l_i$. En effet, p. 132 l.1 de \cite{[T02]}, on peut choisir 
$\Psi _{0}$ non ramifi\'e en les $l_{i}$ (en fait, non ramifi\'e
hors des premiers de $L$ au dessus de $p$). 
Pr\'ecisons le choix du $\Psi_{0}$ l. 13. On ajoute \`a $T$
les places au dessus des $l_{i}$ et on d\'efinit pour ces places
$x$ comme caract\`ere $\Psi_{x}$ le caract\`ere trivial.     
La fin de la preuve du lemme reste la m\^eme.

Dans l'application du th\'eor\`eme de Moret-Bailly p. 136
(th. 1.3. de  \cite{[MB89]} partie 2), on peut imposer que le corps
$E$ est non ramifi\'e en les $l_{i}$.
En effet, il existe une vari\'et\'e ab\'elienne $A$ sur $\Q$ avec 
multplication par $O_{M}$ et de dimension
$[M:\Q]$, principalement polaris\'ee du type
HBAV (p. 133) et qui a bonne r\'eduction en 
les $l_{i}$ : prendre une courbe elliptique
ayant bonne r\'eduction en les $l_{i}$ et tensoriser
par $O_{M}$ (\it cf \rm lemme 1.4. de  \cite{[T02]}).
Les structures de niveau d\'efinissant
le sch\'ema de module $X$ p. 136 sont  non ramifi\'ees en les 
$l_{i}$. La vari\'et\'e ab\'elienne $A$ fournit  
donc un point de $X$ \`a valeurs dans une extension
finie $L_{l_{i}}$ non ramifi\'ee convenable de $\Q _{l_{i}}$.
On en d\'eduit la proposition dans le cas de  \cite{[T02]}.

Passons au cas o\`u l'inertie en $p$ agit sur $\overline{\rho}$
par des caract\`eres de niveau 2 (\cite{[T03]}). 
On choisit :

- p. 29 (de \cite{[T03]}), le corps quadratique
imaginaire $M$ non ramifi\'e en les $l_{i}$ ;

- p. 30, $\chi_{0}$ non ramifi\'e en les $l_{i}$ ; $p_{i}$,
$i=1,2$ distincts des $l_{i}$ ;

- dans le lemme 4.3.,  $\chi$ non ramifi\'e en les $l_{i}$
comme on l'a fait dans le lemme 1.1. de \cite{[T02]}.

Alors, le caract\`ere $\chi_{\lambda}$ d\'efini p. 32 l. -13, -9
(pour $x=\lambda$) est non ramifi\'e en les $l_i$. Les espaces
de modules $X_{\overline{\rho}}$ et $X_{Dih}$ deviennent
isomorphes sur une extension non ramifi\'ee des $\Q _{l_i}$.
On a donc un point de  $X_{\overline{\rho}}$ rationnel
sur un corps $F$ satisfaisant aux conclusions du th\'eor\`eme
de Taylor et qui de plus est non ramifi\'e en les $l_i$.
Cela prouve le compl\'ement.$\Box$

\subsubsection{Le corps $\underline{E}(\rho )$.}\label{Erho}

La proposition suivante est une cons\'equence imm\'ediate
du th\'eor\`eme de Taylor.

\begin{prop}\label{Erhoprop}
Soit $\rho$  comme dans l'\'enonc\'e du th\'eor\`eme
de Taylor. 
Alors, il existe une extension finie $\underline{E}(\rho )$
de $\Q$ contenue dans $E$ caract\'eris\'ee par la propri\'et\'e
suivante :

- il existe une extension finie $L_0$ de $\Q$ contenue
dans $\aQ$ telle que, pour $L$ extension finie de $L_0$
contenue dans $\aQ$, $\underline{E}(\rho )$ soit le sous-corps
de $E$ engendr\'e par les coefficients des polyn\^omes
caract\'eristiques des $\rho (\mathrm{Frob}_ {\mathcal{L}} )$,
pour $\mathcal{L}$ premier de $L$ premier aux \'el\'ements
de $S\cup \{ p \}$.

Pour $F$, $\pi$, $E(\pi )$ et $i$ comme dans l'\'enonc\'e du th\'eor\`eme
de Taylor, $i (E(\pi ))$ contient $\underline{E}(\rho )$.\end{prop} 

\it Preuve de la proposition. \rm
Soient $F$, $\pi$, $E(\pi )$ et  $i$ comme dans l'\'enonc\'e du th\'eor\`eme
de Taylor.  
Alors, pour $L$ extension finie de $F$ contenue dans $\aQ$,
et $\mathcal{L}$ premier de $L$ qui n'est pas au dessus
d'un premier de $S\cup \{ p \}$,
les coefficients des polyn\^omes caract\'eristiques de 
$\rho (\mathrm{Frob}_ {\mathcal{L}} )$ appartiennent
\`a $i( E(\pi ))$. 
Notons  $\underline{E}_{L}$ le sous-corps de 
$i( E(\pi ))$ qu'ils engendrent. 
On pose : $\underline{E}(\rho )=\cap _L \underline{E}_{L}$.
La proposition r\'esulte
imm\'ediatement de ce que, pour $L\subset L'$, on a 
$\underline{E}_{L'}\subset \underline{E}_{L}$.$\Box$

\subsection{Construction d'une repr\'esentation $q$-adique.}\label{qadique}

\begin{prop}\label{propprinc} Supposons $p> 3$. 
Soit $\rho : G_{\Q}\rightarrow \GL _2 (E)$ 
une repr\'esentation $p$-adique absolument
irr\'eductible, g\'eom\'etrique de conducteur $1$
et de poids de Hodge-Tate $(0,k-1 )$. On suppose $k$ pair
(ou, ce qui est \'equivalent, $\rho$ impaire, cf corollaire \ref{det}.
On suppose que $2\leq k \leq (p+1)/2$. Soient $q$ un nombre premier
$\not= 2$,
$\aQ _q$ une cl\^oture alg\'ebrique du corps des nombres $q$-adiques, et
$i_q $ un plongement du corps $\underline{E}(\rho)$ (\ref{Erho})
dans  $\aQ _q$. Alors, il existe une extension 
finie $E_q$ de $\Q _q$ contenant 
$i_q (\underline{E}_{\rho})$ et une repr\'esentation
$q$-adique $\rho_q : G_{\Q}\rightarrow  \GL _2 (E_q)$
qui v\'erifie les deux propri\'et\'es suivantes :

- la restriction de $\rho_q$ \`a tout sous-groupe ouvert de $G_{\Q}$
est absolument irr\'eductible (autrement dit $\rho_{q}$ n'est pas 
potentiellement ab\'elienne)
, $\rho_{q}$ est impaire, 
g\'eom\'etrique de conducteur $1$
et de poids de Hodge-Tate $(0,k-1 )$ ;

- il existe $L$ extension finie de $\Q$ contenue dans 
$\aQ$  telle que pour $\mathcal{L}$ premier
de $L$ qui n'est pas au-dessus de $\{ p,q \}$, les coefficients
$\mathrm{tr}(\rho(\mathrm{Frob}_{\mathcal{L}}))$ et 
$\mathrm{det}(\rho(\mathrm{Frob}_{\mathcal{L}}))$
du  polyn\^ome
caract\'eristique de $\rho (\mathrm{Frob}_{\mathcal{L}})$
appartiennent \`a $\underline{E}(\rho )$ et que l'on ait :

$$\mathrm{tr}(\rho_q(\mathrm{Frob}_{\mathcal{L}}))=
i_q (\mathrm{tr}(\rho(\mathrm{Frob}_{\mathcal{L}}))),\\
\mathrm{det}(\rho_q(\mathrm{Frob}_{\mathcal{L}}))=
i_q (\mathrm{det}(\rho(\mathrm{Frob}_{\mathcal{L}}))).
$$\end{prop}

\it Preuve de la proposition. \rm Soit $\rho$ comme dans l'\'enonc\'e
de la proposition.

\begin{lemme}\label{absirr} Pour toute extension finie $L$ de $\Q$ contenue
dans $\aQ$, la restriction $\rho_{\mid G_L}$ de $\rho $ au groupe de Galois
$G_L$ de $\aQ / L$ est absolument irr\'eductible.\end{lemme}

\it Preuve du lemme. \rm Soit $L$ une extension finie de $\Q$ contenue
dans $\aQ$. La repr\'esentation $\rho_{\mid G_L}$ est semi-simple.
Comme elle est de dimension $2$,
elle est soit  absolument irr\'eductible, soit son image est ab\'elienne.
Comme $\rho$ est absolument irr\'eductible, g\'eom\'etrique de conducteur $1$,
le deuxi\`eme cas
est exclus par la proposition \ref{reppotab}. Le lemme en
r\'esulte. 

\vskip0.5cm

Soit  $\ell$ un nombre premier $\not= p$. Appliquons le th\'eor\`eme
de Taylor (th \ref{taylortheo}) et son compl\'ement
(prop. \ref{complement}) avec $\{ \ell \}$, ce qui est possible puisque
$\rho$ de conducteur $1$ est non ramifi\'ee en $\ell$. 
Soient donc $F^{\{\ell\}}$
un corps totalement r\'eel 
non ramifi\'e au dessus de $\ell$, 
$\pi^{\{\ell\}}$ une repr\'esentation
cuspidale de $\GL _2 (\mathbb{A}_{F^{\{\ell\}}})$,
$i^{\{\ell\}}$ un plongement du corps des coefficients  
$E(\pi^{\{\ell\}}) $ de
$\pi^{\{\ell\}}$ dans $E$ tels que 
la restriction de $\rho$ \`a $G_{F^{\{\ell\}}}$ soit
isomorphe \`a $\rho_{\pi^{\{\ell\}},i^{\{ \ell\}}}$. 

Comme $i^{\{ \ell\}}(E(\pi^{\{\ell\}}))$ contient
les coefficients des polyn\^omes caract\'eristiques
des $\rho (\mathrm{Frob}_{\mathcal{L}} )$ pour 
$\mathcal{L}$ premier de $F^{\{\ell\}}$ qui n'est pas 
au-dessus de $p$,
il r\'esulte de la proposition \ref{Erhoprop} que le corps
$\underline{E}(\rho)$ est contenu dans $i^{\{ \ell\}}(E(\pi^{\{\ell\}}))$.
Choisissons un plongement $\underline{i}^{\{ \ell\}}_q$ de 
$E(\pi^{\{\ell\}})$ dans $\aQ _q$ qui prolonge 
$i_q$, consid\'er\'e  comme un plongement de 
$i^{-1}(\underline{E}(\rho ))$ dans $\aQ _q$ via $i$. 
Notons $E_q ^{\{\ell\}}$
l'adh\'erence de $i^{\{ \ell\}}_{q}(E(\pi^{\{\ell\}}))$
dans $\aQ _q$.
Soit $\rho^{\{\ell\}}_q : G_{F^{\{\ell\}}}\rightarrow \GL _2
( E_q ^{\{\ell\}})$ la repr\'esentation 
associ\'ee \`a $\pi^{\{\ell\}}$ et $\underline{i}^{\{ \ell\}}_q$
par  Taylor
dans \cite{[T89]}. 

Soit $\tau \in G_{\Q}$. Notons  $\mathrm{int}(\tau )(\rho^{\{\ell\}}_q)$
la repr\'esentation 
$q$-adique de $G_{F^{\{\ell\}}}$ conjugu\'ee de $\rho^{\{\ell\}}_q$
par $\tau$.  Dans la repr\'esentation 
$\mathrm{int}({\tau})(\rho^{\{\ell\}}_q)$,
le groupe de Galois $G_{F^{\{\ell\}}}$
agit  donc sur le m\^eme espace 
vectoriel $E_q ^{\{\ell\}}\oplus  E_q ^{\{\ell\}}$ que 
dans la repr\'esentation
$\rho^{\{\ell\}}_q$, mais par 
$\sigma\mapsto \rho^{\{\ell\}}_q (\tau \sigma \tau^{-1})$.

\begin{lemme}\label{conjegal} La repr\'esentation 
$\mathrm{int}(\tau )(\rho^{\{\ell\}}_q)$
est isomorphe \`a $\rho^{\{\ell\}}_q$ : il existe $g_{\tau}\in
\GL _2 (E_q ^{\{\ell\}})$ tel que, pour tout $\sigma\in G_{F^{\{\ell\}}}$ :

$$ \rho^{\{\ell\}}_q (\tau \sigma \tau ^{-1})=
g_{\tau}\rho^{\{\ell\}}_q(\sigma )g_{\tau}^{-1}.$$

\end{lemme}

\it Preuve. \rm
La  repr\'esentation galoisienne $\mathrm{int}(\tau )(\rho^{\{\ell\}}_q)$
est associ\'ee \`a la repr\'esentation
automorphe $^{\tau}\pi^{\{\ell\}}$
obtenue en composant $\pi^{\{\ell\}}$ avec 
l'automorphisme de $\GL_2 (\mathbb{A}_{F^{\{\ell\}}})$ d\'efini par $\tau$
et au plongement $i$ du corps de d\'efinition 
$E( ^{\tau}\pi^{\{\ell\}} )$
de $^{\tau}\pi^{\{\ell\}}$
(clairement $E( ^{\tau}\pi^{\{\ell\}} )= E( \pi^{\{\ell\}} )$).
La repr\'esentation $p$-adique associ\'ee \`a 
$^{\tau}\pi^{\{\ell\}}$ est la repr\'esentation 
conjugu\'ee $\mathrm{int}( \tau )(\rho_{F^{\{\ell\}}})$ de la restriction
$\rho_{F^{\{\ell\}}}$
de $\rho$ \`a $G_{F^{\{\ell\}}}$. La repr\'esentation  
$\mathrm{int}(\tau)(\rho_{F^{\{\ell\}}})$ est isomorphe \`a 
$\rho_{F^{\{\ell\}}}$. En effet, pour tout $\sigma\in G_{F^{\{\ell\}}}$  :

$$\mathrm{int}(\tau )(\rho_{F^{\{\ell\}}}(\sigma))=\rho (\tau) 
\rho_{F^{\{\ell\}}}(\sigma)\rho (\tau)^{-1}.$$

Il r\'esulte alors du th\'eor\`eme de multiplicit\'e $1$ fort
(\cite{[PS81]})
que $^{\tau}\pi^{\{\ell\}}$ est isomorphe \`a $\pi^{\{\ell\}}$.
On en d\'eduit le lemme.

\begin{lemme}\label{rhoqgl2} Pour toute extension finie $L$
de $F^{\{\ell\}}$ contenue dans $\aQ$, la restriction de $\rho^{\{\ell\}}_q$
\`a $G_L$ est absolument irr\'eductible.
\end{lemme}

\it Preuve. \rm 
Soit $L$ une extension finie de $F^{\{\ell\}}$ contenue dans $\aQ$.
Soit $(\rho^{\{\ell\}}_q)_{\mid G_L}$ est  absolument irr\'eductible,
soit sa semi-simplifi\'ee a une image ab\'elienne.
Le second cas est exclus. En effet, les repr\'esentations $p$-adiques 
et $q$-adiques  $\rho_{\mid G_{L}}$        et $(\rho^{\{\ell\}}_q)_{\mid G_L}$ 
sont compatibles (1.2.3. de \cite{[SAB]}).
Si $(\rho^{\{\ell\}}_q)_{\mid G_L}$ a une image ab\'elienne,
il r\'esulte du corollaire 1 du th. 2 du 2.3. de 
\cite{[SAB]} que la restriction de 
$\rho$ \`a $G_{L}$ a aussi une image ab\'elienne.
Ce n'est pas le cas (proposition \ref{potab}). Le lemme
est prouv\'e.

\vskip0.5cm

En particulier, la rep\'esentation $\rho^{\{\ell\}}_q$ est  absolument 
irr\'eductible. 
Si $g_{\tau}$ est comme dans le lemme \ref{conjegal},
il en r\'esulte que 
l'image $\overline{g_{\tau}}$ de $g_{\tau}$ dans 
$\mathrm{PGL}_2 (E_q ^{\{\ell\}} )$ ne d\'epend pas du choix 
de $g_{\tau}$. On voit alors que $\tau\mapsto 
\overline{g_{\tau}}$ est une repr\'esentation projective 
de $G_{\Q}$. On la note $\rho_{q,\mathrm{proj}}^{\{\ell\}}$.

\begin{lemme}\label{indell} La restriction de  
$\rho_{q,\mathrm{proj}}^{\{\ell\}}$ \`a $G_{F^{\{\ell\}}}$
est le compos\'e de $\rho_{q}^{\{\ell\}}$ avec 
la projection $\GL _2 (E_q ^{\{\ell\}} )\rightarrow 
\mathrm{PGL}_2 (E _q ^{\{\ell\}} )$. La repr\'esentation 
$\rho_{q,\mathrm{proj}}^{\{\ell\}}$ ne d\'epend pas 
du choix de $\ell$ : plus pr\'ecis\'ement, si 
$\ell_1$ est un premier $\not= p$, et si $E _q ^{\{\ell, \ell_1 \}}$
est le compos\'e dans $\aQ _q$ des corps $E_q ^{\{\ell\}}$ et 
$E_q ^{\{\ell _1\}}$, il existe un unique 
$g\in \mathrm{PGL}_2 (E _q ^{\{\ell , \ell_1 \}})$ tel que 
$\rho_{q,\mathrm{proj}}^{\{\ell _1 \}}=\mathrm{int}(g )
(\rho_{q,\mathrm{proj}}^{\{\ell  \}})$.
\end{lemme}

\it Preuve. \rm La premi\`ere partie du lemme est claire. 

Prouvons la seconde partie.
Soit $\ell _1$ un premier premiers $\not= p$.
Soit $L$ une extension galoisienne finie de $\Q$ contenue dans $\aQ$
contenant $F^{\{\ell\}}$ et $F^{\{\ell _1\}}$, et telle
que les coefficients des polyn\^omes caract\'eristiques
des $\rho (\mathrm{Frob}_{\mathcal{L}} )$, pour 
$\mathcal{L}$ premier de $L$ premier \`a $p$, soit
contenus dans $\underline{E}(\rho )$.   
Les restrictions de $\rho_{q}^{\{\ell \}}$ et 
$\rho_{q}^{\{\ell_1 \}}$ \`a $G_L$
sont absolument irr\'eductibles (lemme pr\'ec\'edent).
Elles ont m\^eme caract\`ere puisque, pour $\mathcal{L}$
premier de $L$ premier \`a $\ell$, $\ell _1$ et $p$ : 

$$\mathrm{tr}(\rho_{q}^{\{\ell \}} (\mathrm{Frob}_{\mathcal{L}} ))=
\mathrm{tr}(\rho_{q}^{\{\ell _1 \}} (\mathrm{Frob}_{\mathcal{L}} ))=
i_q (\mathrm{tr}(\rho (\mathrm{Frob}_{\mathcal{L}} ))).$$

Elles sont donc isomorphes et  il existe un unique 
$g\in  \mathrm{PGL}_2 (E _q ^{\{\ell , \ell_1 \}})$ tel que :

$$ (\rho_{q}^{\{\ell_1 \}})_{\mid G_L}=\mathrm{int}(g) 
(\rho_{q}^{\{\ell \}})_{\mid G_L}.$$

Notons $\eta$ cette repr\'esentation de $G_L$.
Soit $\tau \in G_{\Q}$. Posons 
$\overline{g_0 }=
\rho_{q,\mathrm{proj}}^{\{\ell _{1} \}}(\tau )$ et 
$\overline{g_1 }=
\mathrm{int}(g)(\rho_{q,\mathrm{proj}}^{\{\ell \}}(\tau ))$.
Pour $i=0,1$, $\overline{g_i }$
v\'erifie, pour tout $\sigma \in G_L$ :

$$ \mathrm{int}(\overline{g_i })(\eta (\sigma ))
=\eta (\mathrm{int}(\tau )(\sigma )).$$  

Comme ceci caract\'erise $\overline{g_i }$, on a  
$\overline{g_0 }=\overline{g_1 }$, ce qui prouve le lemme.

\vskip0.5cm

On choisit $\ell\not=p$ et on  
on pose $\rho_{q,\mathrm{proj}}
=\rho_{q,\mathrm{proj}}^{\{\ell\}}$. 
On pose $E_q =E_q^{\{ \ell \}}$.
Donc  $\rho_{q,\mathrm{proj}}$
est \`a valeurs  dans $\mathrm{PGL}_2 (E_q )$.
Pour tout $\ell _1\not=p$, on identifie 
$\rho_{q,\mathrm{proj}}^{\{\ell\}}$ \`a $\rho_{q,\mathrm{proj}}$
gr\^ace au lemme pr\'ec\'edent.

\begin{lemme}\label{ramif} La repr\'esentation 
$\rho_{q,\mathrm{proj}}$
est non ramifi\'ee en dehors de $q$. Soit $I_q \subset G_{\Q}$
un sous-groupe d'inertie en $q$ ; comme $F^{\{q\}}$ est non
ramifi\'e en $q$,  $I_q$ s'identifie \`a un sous-groupe 
de $G_{F^{\{q\}}}$. La restriction de 
$\rho_{q,\mathrm{proj}}$ \`a $I_q$ co\"{\i}ncide 
\`a un automorphisme int\'erieur pr\`es
avec 
le compos\'e de $(\rho_{q}^{\{q\}})_{\mid I_q}$ avec la projection 
$\GL _2 (E_q )\rightarrow \mathrm{PGL}_2 (E_q )$. Si $q\not= 2$, 
la repr\'esentation $(\rho_{q}^{\{q\}})_{\mid I_q}$ 
est cristalline de type de Hodge-Tate 
constant  $(0,k-1)$ (prop. \ref{ht}).
\end{lemme}

\it Preuve. \rm 
Si $q=p$, la proposition est claire.  On suppose $q\not= p$.

Soit $\ell _1$ un premier
$\not =p$.  Comme $\rho$ est non ramifi\'ee en dehors de $p$,
il r\'esulte d'un th\'eor\`eme de H Carayol (\cite{[C86]})
compl\'et\'e par R Taylor (introduction de \cite{[T89]})
que la repr\'esentation automorphe $\pi ^{\{\ell _1\}}$ est
non ramifi\'ee en $\ell _1$. Si $\ell _1 \not=q$, on 
voit que  $\rho^{\{\ell _1\}}_q$ est 
non ramifi\'ee en $\ell _1$. Comme $F^{\{\ell _1\}}$
est non ramifi\'e en $\ell _1$,  $\rho_{q,\mathrm{proj}}$
est non ramifi\'ee en $\ell _1$ ($\ell _1\not=p,q$).
Elle est aussi non ramifi\'ee en $p$ 
car $\pi ^{\{\ell\}}$ et $F^{\{\ell\}}$ le sont
(th. \ref{taylortheo}) et  que l'on a suppos\'e 
$q\not=p $. On a donc prouv\'e que $\rho_{q,\mathrm{proj}}$
est non ramifi\'ee en dehors de $q$.

Comme $\rho$ est non ramifi\'ee en $q$,
$\pi^{\{q\}}$ est non ramifi\'ee en $q$.
Comme $F^{\{q\}}$ est non ramifi\'e en $q$ et que $q\not= 2$, 
il r\'esulte de C Breuil (\cite{[BR99]}) que la restriction de
$\rho^{\{q\}}_q$ \`a un sous-groupe d'inertie en 
un premier de $F^{\{q\}}$ au dessus de $q$  
est limite $q$-adique de repr\'esentations 
cristallines de poids  de Hodge-Tate  $\in \{ 0,k-1 \} $.
Il r\'esulte de L Berger (\cite{[Ber]}) que  $\rho^{\{q\}}_q$
est cristalline de poids de Hodge-Tate $\in \{ 0,k-1 \} $.
Comme  le caract\`ere de Dirichlet de la forme
modulaire de Hilbert associ\'ee \`a $\pi^{\{q\}}$ est trivial,
le d\'eterminant de $\rho^{\{q\}}_q$ est la restriction
de $ \chi_q ^{k-1}$ \`a $G_{F^{\{q\}}}$. On voit alors que 
le type de $\rho^{\{q\}}_q$ est bien constant \'egal
\`a $(0,k-1)$. Ceci ach\`eve la preuve du lemme.$\Box$

\vskip0.5cm

Notons $O_q$ l'anneau de valuation de $E_q$.
L'image de $\rho_{q,\mathrm{proj}}$ est compacte. Les sous-groupes compacts 
maximaux de $\mathrm{PGL}_2 (E_q )$
sont les images des sous-groupes du type  $\GL (T)$,
pour $T\subset (E_q )^2$ un $O_q$-r\'eseau. Apr\`es conjugaison,
on peut supposer que l'image de $\rho_{q,\mathrm{proj}}$ 
est contenue dans $\mathrm{PGL}_2 (O_q )$. 

\begin{lemme} 
La repr\'esentation $\rho_{q,\mathrm{proj}} :
G_{\Q}\rightarrow \mathrm{PGL}_{2} (O_{q})$
se rel\`eve en une repr\'esentation $\rho_q: G_{\Q}
\rightarrow \GL _2 (E_q )$ qui est irr\'eductible, 
non ramifi\'ee en dehors de $q$, et si $q\not= 2$ cristalline en $q$
de poids de Hodge-Tate $(0,k-1)$.
\end{lemme}

\vskip0.5cm 

\it Preuve du lemme. \rm 
Notons $d$ le morphisme $\mathrm{PGL}_2 (O_q )
\rightarrow O_q ^* /(O_q ^*)^2$ d\'efini par le d\'eterminant.
Notons $G$ le sous-groupe de 
$\mathrm{PGL}_2 (O_q )\times O_q ^*$ form\'e des 
$(g,\lambda)$ tels que 
$d(g)=\lambda\ \mathrm{mod}(O_q ^*)^2$. On a la suite exacte :

$$ (1)\quad : \quad 1 \rightarrow \{ \pm 1 \}\rightarrow \GL _2 (O_q )
\rightarrow G \rightarrow 1.$$

Notons $k_{q}$ le corps r\'esiduel de $O_{q}$.
De m\^eme, le d\'eterminant d\'efinit un morphisme $\overline{d} :
\mathrm{PGL}_2 (k_q )\rightarrow 
k_q ^* /(k_{q} ^*)^2$. Notons
 $\overline{G}$ le sous-groupe
de $\mathrm{PGL}_2 (k_q )\times k_q ^*$
form\'e des $(\overline{g},\overline{\lambda})$
tels que $\overline{d}(\overline{g})=
\overline{\lambda}\ \mathrm{mod}(\mathbb{F}_q ^*)^2$.
On a la suite exacte :

$$ (\overline{1})\quad : \quad 1 \rightarrow \{ \pm 1 \}\rightarrow \GL _2 
(k_q )
\rightarrow \overline{G} \rightarrow 1.$$

Comme $q\not= 2$, la suite exacte $(\overline{1})$ est le ``pull-back'' de la 
suite $(1)$.
La repr\'esentation $\gamma$ :

$$\rho_{q,\mathrm{proj}} \times \chi _q^{k-1}: G_{\Q}
\rightarrow \mathrm{PGL}_2 (O_q )\times O_q ^*$$

est \`a valeurs dans $G$. 
Soit en effet $\delta: G_{\Q}\rightarrow k_q ^* /
(k_q ^*)^2$ le compos\'e de $\gamma$ avec
$(\overline{g},\lambda)\mapsto \overline{d(g)}
\overline{\lambda}^{-1}$. Comme $\gamma$ est non ramifi\'e
en dehors de $q$, le caract\`ere $\delta$ est non ramifi\'e
en dehors de $q$. Il est aussi non ramifi\'e en $q$.
En effet, le lemme pr\'ec\'edent entra\^{\i}ne que la 
restriction de $\delta$ \`a $I_q$ est triviale.
On voit  que $\delta$ est non ramifi\'e partout,
donc trivial et $\gamma$ est bien \`a valeurs dans 
$G$.

Par r\'eduction modulo $q$, on obtient une repr\'esentation
$\overline{\gamma}$ de $G_{\Q}$ dans $\overline{G}$.
L'obstruction \`a relever $\overline{\gamma}$ en une repr\'esentation
$G_{\Q }\rightarrow \GL _2 (k_q )$ est un 
\'el\'ement $o(\overline{\gamma})$ dans 
 le groupe $\mathrm{Br} (\Q )_2$ des \'el\'ements tu\'es par $2$
du groupe de Brauer $\mathrm{Br} (\Q )$. 

Cette obstruction
v\'erifie le principe de Hasse. Pour chaque premier $\ell$,
la composante  locale $o(\overline{\gamma})_{\ell}$ est 
l'obstruction \`a relever la restriction de $\overline{\gamma}$
\`a un sous-groupe de d\'ecomposition en $\ell$. Si $\ell \not= q$,
elle est nulle car $\overline{\gamma}$ est non ramifi\'ee en $\ell$.

L'obstruction $o(\overline{\gamma})_q$ est aussi nulle.
Soit en effet $I_q \subset G_{\Q}$ un sous-groupe d'inertie
en $q$. La restriction de $\rho^{\{ q\}}_q : G_{F^{\{ q\}}}
\rightarrow \GL _2 (E_q )$ \`a $I_q$ est un rel\`evement 
de la restriction de $\gamma$ \`a $I_q$. Il est cristallin.
L'unicit\'e de ce rel\`evement cristallin 
entra\^{\i}ne alors qu'il
se prolonge au sous-groupe de d\'ecomposition en $q$
(proposition 1.2. de \cite{[W95]}).

On a donc que les composantes locales $o(\overline{\gamma})_{\ell}$,
pour $\ell$ premier de $\Q$, sont nulles. C'est donc que 
$o(\overline{\gamma})$ est nulle. La repr\'esentation 
$\overline{\gamma }$ se rel\`eve en une repr\'esentation
$G_{\Q}\rightarrow \GL _2 (\mathbb{F}_q )$. Il en r\'esulte
que $\gamma$ se rel\`eve en une repr\'esentation $\widehat{\gamma} :
G_{\Q}\rightarrow \GL _2 (O_q)$ qui est non ramifi\'ee
en dehors d'un ensemble fini $S$ de nombres premiers.
De plus, pour chaque premier $\ell\in S$, si $I_{\ell}\subset G_{\Q}$ est 
un sous-groupe d'inertie en $\ell$, il existe 
un caract\`ere $\eta_{\ell}: I_{\ell}\rightarrow \{ \pm 1\}$
qui est tel que $\widehat{\gamma}_{\mid I_{\ell}} \eta_{\ell}$ soit 
non ramifi\'ee si $\ell \not= q$, et cristalline si $\ell =q$.
Pour $\ell\not =2$, soit $L_{\ell}=\Q (\sqrt{\epsilon_{\ell }\ell})$, 
$\epsilon_{\ell}=(-1)^{l-1/2}$. Le corps quadratique $L_{\ell}$
 est ramifi\'ee seulement en $\ell$ et
$I_{\ell}\rightarrow \mathrm{Gal}(L_{\ell}/\Q)$ 
identifie $\mathrm{Gal}(L_{\ell}/\Q)$ avec le plus grand 
quotient de $I_{\ell}$ qui est un groupe ab\'elien 
tu\'e par $2$.
Soit de m\^eme, $L_2 =\Q (\sqrt{2},i)$ ; $L_2$  est non ramifi\'e 
en dehors de $2$ et
$I_2\rightarrow  
\mathrm{Gal}(L_2 / \Q)$  identifie $\mathrm{Gal}(L_2 / \Q)$
avec le plus grand quotient de $I_2$ qui est un groupe
ab\'elien tu\'e par $2$. On voit alors qu'il existe 
un caract\`ere $\eta : G_{\Q}\rightarrow \{ \pm 1\}$ 
qui est tel que $\widehat{\gamma}\eta$
soit non ramifi\'e en dehors de $q$ et cristalline
en $q$. 

Posons $\rho_{q}= \widehat{\gamma}\eta$ ; $\rho_{q}$ est donc 
g\'eom\'etrique de conducteur $1$. 
Soit $\ell$ un nombre premier $\not= p$. 
Il existe un caract\`ere \`a image finie $\eta$ de $G_{F^{\{ \ell \}}}$
tel que la restriction
de $\rho_{q}$ \`a $G_{F^{\{ \ell \}}}$ co\"{\i}ncide
avec la repr\'esentation obtenue \`a partir de  
$\rho_{q}^{\{ \ell \}}$ par torsion par le  caract\`ere $\eta$.
Il en r\'esulte que $\rho_{q}$ est impaire. 
Il r\'esulte aussi du lemme \ref{rhoqgl2} que la restriction
de $\rho_{q}$ \`a tout sous-groupe ouvert de $G_{\Q}$ est 
absolument irr\'eductible. Enfin, la propri\'et\'e de compatibilit\'e
avec $\rho$ de la proposition est v\'erifi\'ee pour $L$ correspondant
au noyau de $\eta$. La proposition est prouv\'ee.

\subsection{Conclusion de la
d\'emonstration du th\'eor\`eme \ref{theoprinc}.}.\label{conclusion}

Soit donc $\rho$ une repr\'esentation $p$-adique g\'eom\'etrique
de conducteur $1$ et de poids $(0,k)$, $k=2$, ou
$k=4$ et $p\geq 7$.

Supposons tout d'abord $p>3$.
Supposons $\rho$ irr\'eductible. Il r\'esulte de la proposition
\ref{reppotab} que $\rho$ n'est pas potentiellement ab\'elienne.
On applique alors la proposition \ref{propprinc}
avec $q=7$.  
JM Fontaine a  prouv\'e 
qu'une repr\'esentation $p$-adique 
$\rho$ de $G_{\Q}$ dans un $\Qp$-espace vectoriel de dimension
finie arbitraire, qui est non ramifi\'ee en dehors de $p$
et cristalline en $p$ 
\`a poids de Hodge-Tate dans l'intervalle $[0,h]$
est extensions de $\Qp (i)$  
sous l'une ou l'autre des deux hypoth\`eses suivantes  :

- $h=1$ et $2<p\leq 17$ (\cite{[F85]});

- $h=3$ et $p=7$ (\cite{[F93]}). 

V Abrashkin a donn\'e ind\'ependamment de JM Fontaine 
une d\'emonstration du cas $k=2$ qui marche aussi 
pour $p=2$ (\cite{[A89]}).
On voit donc que $\rho_{7}$ en fait n'existe pas. Donc $\rho$
n'est pas irr\'eductible. 

Si $k=2$ et $p=2,3$, les th\'eor\`emes de Fontaine et Abrashkin
donnent directement que $\rho$ est r\'eductible. 

On voit donc que $\rho$  est r\'eductible. Les r\'esultats de 
Soul\'e cit\'es dans l'introduction entra\^{\i}nent alors le 
th\'eor\`eme.

\section{Application aux groupes $p$-divisibles.}\label{grpdiv}

\begin{cor} Soit $\Gamma$ un groupe $p$-divisible sur
$\Z$. Soient 
$O_E$ l'anneau des entiers de $E$ et $O\subset O_E$
un ordre de $O_E$. Soit $n$ le degr\'e de $E/\Qp$.
On suppose que $\Gamma$ est de hauteur $2n$ et que
l'on a un plongement de $O$ dans $\mathrm{End} (\Gamma )$.
Alors, $\Gamma$ est isog\`ene \`a une somme directe
de deux groupes $p$-divisibles $\Gamma _1$ et 
$\Gamma _2$, $\Gamma _1$ et 
$\Gamma _2$ \'etant isomorphes soit au groupe constant 
$E/ O_E$ ou \`a son dual de Cartier $O_E \otimes \mu_{p^{\infty}}$.
\end{cor}

\it D\'emonstration du corollaire. \rm 
Soient $T$ le module de Tate de $\Gamma$ et 
$U=\Qp \otimes_{\Zp} T$. La repr\'esentation $\rho (\Gamma )$ de
$G_{\Q}$ dans $U$ est non ramifi\'ee 
en dehors de $p$ et cristalline de poids $(0,0)$,
$(0,1)$ ou $(1,1)$. Si les poids sont $(0,0)$,
$\rho (\Gamma )$ est triviale. Si elle est de poids
$(1,1)$, la repr\'esentation tordue  $U(-1)$ l'est. Si elle est de poids 
$0$ et $1$, le th\'eor\`eme entra\^{\i}ne qu'elle
est isomorphe \`a $E\oplus E(1)$. Ceci  prouve le corollaire,
puisque $\rho (\Gamma )$ d\'etermine $\Gamma$ \`a isog\'enie pr\`es
(\cite{[Ta67]}).

\section{Repr\'esentations g\'eom\'etriques de conducteur $1$
de bas poids ordinaires en $3$.}\label{ord}

\begin{prop} Soient $p$ un premier $>2$ et $k$ un entier
v\'erifiant $2\leq k \leq (p+1)/2$. Il n'existe pas de 
repr\'esentations $\rho : G_{\Q} \rightarrow \GL _{2} (E)$
qui soit impaire, irr\'eductible, g\'eom\'etrique de conducteur
$1$ et de poids $k$, et qui v\'erifie l'hypoth\`ese d'ordinarit\'e
en $3$ suivante :

- il existe un plongement $i_{3}: \underline{E} (\rho )
\hookrightarrow \aQ _{3}$ et une valeur propre $\lambda$
de $\rho (\mathrm{Frob}_{3} )$ telle que 
$i_{3}\circ i^{-1} (\lambda )$ soit une unit\'e $3$-adique
(pour la d\'efinition de $\underline{E} (\rho )$,
voir \ref{Erho}).\end{prop}

\it Preuve. \rm Soit $\rho$ comme dans l'\'enonc\'e de la 
proposition. 
On applique la proposition \ref{propprinc} avec 
$q=3$ et $i_{3}$:
on obtient la repr\'esentation $3$-adique $\rho _{3}$.
Soit $\overline{\rho}_{3}$ une r\'eduction de $\rho _{3}$.
Elle  est non ramifi\'ee hors de $3$. Elle est impaire.
Un th\'eor\`eme de Serre dit alors que sa semi-simplifi\'ee 
est isomorphe \`a $1\oplus \overline{\chi _{3}}$ 
(p. 710 de \cite{[Se]}). Un th\'eor\`eme
de Wiles dit que l'hypoth\`ese d'ordinarit\'e en $3$ entra\^{\i}ne que 
la restriction de $\rho_{3}$ au groupe d'inertie en $3$ est du type :
$$\left( \begin{array}{cc}
  	\chi_{3}^{k-1} & *  \\
  	0 & 1
  \end{array} \right),$$

(th. 2 de \cite{[W88]}). 
Un th\'eor\`eme de Skinner et Wiles donne  la modularit\'e de 
$\rho _{3}$ (\cite{[SW99]}). La th\'eorie de Hida dit alors que $\rho_{3}$
n'existe pas ( 7.6. de \cite{[Hida]}). La repr\'esentation  $\rho_{3}$ 
proviendrait d'une forme parabolique $f$
de niveau $1$ ordinaire en $3$, qui serait une sp\'ecialisation
d'une famille de Hida, dont la sp\'ecialisation en poids $2$
serait une forme parabolique de poids $2$ pour $\Gamma _{0} (3)$.
Une telle forme n'existe pas. 

\it Remarque. \rm Si l'on  connaissait le th\'eor\`eme de Skinner et
Wiles sans l'hypoth\`ese d'ordinarit\'e de l'action de $I_{3}$
sur la repr\'esentation $3$-adique, on aurait la modularit\'e
de $\rho$ sans l'hypoth\`ese d'ordinarit\'e en $3$ de $\rho$.

\nocite{*}
\bibliographystyle{plain}
\bibliography{grpepdiv4}

\begin{thebibliography}{10}

\bibitem{[Per]}
{\em P\'eriodes {$p$}-adiques}.
\newblock Soci\'et\'e Math\'ematique de France, Paris, 1994.
\newblock Papers from the seminar held in Bures-sur-Yvette, 1988, Ast\'erisque
  No. 223 (1994).

\bibitem{[A89]}
V.~A. Abrashkin.
\newblock Galois modules of group schemes of period {$p$} over the ring of
  {W}itt vectors.
\newblock {\em Izv. Akad. Nauk SSSR Ser. Mat.}, 51(4):691--736, 910, 1987.

\bibitem{[BN]}
Denis Benois and Thong Nguyen Quang~Do.
\newblock Les nombres de {T}amagawa locaux et la conjecture de {B}loch et
  {K}ato pour les motifs {$\Bbb Q(m)$} sur un corps ab\'elien.
\newblock {\em Ann. Sci. \'Ecole Norm. Sup. (4)}, 35(5):641--672, 2002.

\bibitem{[Ber]}
Laurent Berger.
\newblock Limites de repr\'esentations cristallines.
\newblock {\em A para\^{\i}tre \`a Compositio}.

\bibitem{[BR99]}
Christophe Breuil.
\newblock Une remarque sur les repr\'esentations locales {$p$}-adiques et les
  congruences entre formes modulaires de {H}ilbert.
\newblock {\em Bull. Soc. Math. France}, 127(3):459--472, 1999.

\bibitem{[C86]}
Henri Carayol.
\newblock Sur les repr\'esentations {$l$}-adiques associ\'ees aux formes
  modulaires de {H}ilbert.
\newblock {\em Ann. Sci. \'Ecole Norm. Sup. (4)}, 19(3):409--468, 1986.

\bibitem{[Dieu]}
Luis Dieulefait.
\newblock Existence of families of {G}alois representations and new cases of
  the {F}ontaine-{M}azur conjecture.
\newblock {\em arXiv:math.NT0304433}, 2003.

\bibitem{[F78]}
Jean-Marc Fontaine.
\newblock Modules galoisiens, modules filtr\'es et anneaux de
  {B}arsotti-{T}ate.
\newblock In {\em Journ\'ees de G\'eom\'etrie Alg\'ebrique de Rennes. (Rennes,
  1978), Vol. III}, volume~65 of {\em Ast\'erisque}, pages 3--80. Soc. Math.
  France, Paris, 1979.

\bibitem{[F85]}
Jean-Marc Fontaine.
\newblock Il n'y a pas de vari\'et\'e ab\'elienne sur {${\bf Z}$}.
\newblock {\em Invent. Math.}, 81(3):515--538, 1985.

\bibitem{[F93]}
Jean-Marc Fontaine.
\newblock Sch\'emas propres et lisses sur {${\bf Z}$}.
\newblock In {\em Proceedings of the Indo-French Conference on Geometry
  (Bombay, 1989)}, pages 43--56, Delhi, 1993. Hindustan Book Agency.

\bibitem{[FM95]}
Jean-Marc Fontaine and Barry Mazur.
\newblock Geometric {G}alois representations.
\newblock In {\em Elliptic curves, modular forms, \& Fermat's last theorem
  (Hong Kong, 1993)}, Ser. Number Theory, I, pages 41--78. Internat. Press,
  Cambridge, MA, 1995.

\bibitem{[Hida]}
Haruzo Hida.
\newblock {\em Elementary theory of {$L$}-functions and {E}isenstein series},
  volume~26 of {\em London Mathematical Society Student Texts}.
\newblock Cambridge University Press, Cambridge, 1993.

\bibitem{[MB89]}
Laurent Moret-Bailly.
\newblock Groupes de {P}icard et probl\`emes de {S}kolem. {I}, {II}.
\newblock {\em Ann. Sci. \'Ecole Norm. Sup. (4)}, 22(2):161--179, 181--194,
  1989.

\bibitem{[BPR]}
Bernadette Perrin-Riou.
\newblock La fonction {$L$} {$p$}-adique de {K}ubota-{L}eopoldt.
\newblock In {\em Arithmetic geometry (Tempe, AZ, 1993)}, volume 174 of {\em
  Contemp. Math.}, pages 65--93. Amer. Math. Soc., Providence, RI, 1994.

\bibitem{[PS81]}
I.~I. Piatetski-Shapiro.
\newblock Multiplicity one theorems.
\newblock In {\em Automorphic forms, representations and $L$-functions (Proc.
  Sympos. Pure Math., Oregon State Univ., Corvallis, Ore., 1977), Part 1},
  Proc. Sympos. Pure Math., XXXIII, pages 209--212. Amer. Math. Soc.,
  Providence, R.I., 1979.

\bibitem{[Se]}
Jean-Pierre Serre.
\newblock {\em {\OE}uvres. {V}ol. {III}}.
\newblock Springer-Verlag, Berlin, 1986.
\newblock 1972--1984.

\bibitem{[SAB]}
Jean-Pierre Serre.
\newblock {\em Abelian {$l$}-adic representations and elliptic curves},
  volume~7 of {\em Research Notes in Mathematics}.
\newblock A K Peters Ltd., Wellesley, MA, 1998.
\newblock With the collaboration of Willem Kuyk and John Labute, Revised
  reprint of the 1968 original.

\bibitem{[SW99]}
C.~M. Skinner and A.~J. Wiles.
\newblock Residually reducible representations and modular forms.
\newblock {\em Inst. Hautes \'Etudes Sci. Publ. Math.}, (89):5--126 (2000),
  1999.

\bibitem{[SW01]}
C.~M. Skinner and Andrew~J. Wiles.
\newblock Nearly ordinary deformations of irreducible residual representations.
\newblock {\em Ann. Fac. Sci. Toulouse Math. (6)}, 10(1):185--215, 2001.

\bibitem{[Sk03]}
Chris Skinner.
\newblock Modularity of {G}alois representations.
\newblock {\em J. Th\'eor. Nombres Bordeaux}, 15(1):367--381, 2003.
\newblock Les XXII\`emes Journ\'ees Arithmetiques (Lille, 2001).

\bibitem{[Ta67]}
J.~T. Tate.
\newblock {$p-divisible$} {$groups.$}.
\newblock In {\em Proc. Conf. Local Fields (Driebergen, 1966)}, pages 158--183.
  Springer, Berlin, 1967.

\bibitem{[T89]}
Richard Taylor.
\newblock On {G}alois representations associated to {H}ilbert modular forms.
\newblock {\em Invent. Math.}, 98(2):265--280, 1989.

\bibitem{[T03]}
Richard Taylor.
\newblock On the meromorphic continuation of degree two {L}-functions.
\newblock {\em Preprint}, pages 1--53, 2001.

\bibitem{[T02]}
Richard Taylor.
\newblock Remarks on a conjecture of {F}ontaine and {M}azur.
\newblock {\em J. Inst. Math. Jussieu}, 1(1):125--143, 2002.

\bibitem{[W88]}
A.~Wiles.
\newblock On ordinary {$\lambda$}-adic representations associated to modular
  forms.
\newblock {\em Invent. Math.}, 94(3):529--573, 1988.

\bibitem{[W95]}
J.-P. Wintenberger.
\newblock Rel\`evement selon une isog\'enie de syst\`emes {$l$}-adiques de
  repr\'esentations galoisiennes associ\'es aux motifs.
\newblock {\em Invent. Math.}, 120(2):215--240, 1995.

\end{thebibliography}

\vspace{2cm}
Jean-Pierre Wintenberger

Universit\'e Louis Pasteur

D\'epartement de Math\'ematiques, IRMA

7, rue Ren\'e Descartes 

67084 Strasbourg Cedex

France

e-mail wintenb@math.u-strasbg.fr

tel 03 90 24 02 17 , fax 03 90 24 03 28

\end{document}